\documentclass{amsart}
\usepackage{amsthm,amssymb,mathtools}
\usepackage{stmaryrd,tikz,pgfplots}
\usepackage{marginnote}
\usepackage{enumitem}

\usepackage[textsize=tiny]{todonotes}
\newcommand{\Marginpar}[1]{\marginpar{\tiny{#1}}}
\newcommand{\Note}[1]{{\par\noindent\hrulefill\par\tiny{#1}\par\noindent\hrulefill\par}}
\newcommand{\Detail}[1]{{#1}}
\renewcommand{\Marginpar}[1]{}
\renewcommand{\Note}[1]{}
\renewcommand{\Detail}[1]{}

\usepackage{hyperref}
\usepackage[all]{xy}
\usepackage{amsmath,amsfonts}	% Bunch of math support
\usepackage{amsthm,latexsym,amssymb}
\usepackage{tensor}

\newtheorem{thm}{Theorem}[section]

\theoremstyle{definition}

\renewcommand{\[}{\begin{equation*}}
\renewcommand{\]}{\end{equation*}}

\begin{document}
\parskip1mm

\title[Einstein--Maxwell equations]{Einstein--Maxwell equations on $4$-dimensional Lie algebras}

\author{Caner Koca}
\address{Department of Mathematics, NYC College of Technology of CUNY, Brooklyn, NY 11021, USA.}
\email{ckoca@citytech.cuny.edu}

\author{Mehdi Lejmi}
\address{Department of Mathematics, Bronx Community College of CUNY, Bronx, NY 10453, USA.}
\email{mehdi.lejmi@bcc.cuny.edu}

\begin{abstract}
We classify up to automorphisms all left-invariant non-Einstein solutions to the Einstein--Maxwell equations on $4$-dimensional
Lie algebras.
\end{abstract}
\maketitle

%\setcounter{tocdepth}{1}
%\tableofcontents

\section{Introduction}
A Riemannian $4$-manifold $(M,g)$ is called \emph{Einstein} if the trace-free Ricci tensor is identically zero, that is, $Ric_0 :=Ric - \frac{s}{4} g = 0$. From General Relativity's viewpoint, these are the Riemannian solutions of Einstein's Field Equations in vacuum. One can also consider the same equations in the presence of an \emph{electro-magnetic field} $F$. In Physics, $F$ can be thought as a differential 2-form, which is closed and co-closed: $dF=0$ and $d\star F=0$, where $\star$ is the Hodge-star operator (in particular, the manifold is assumed to be oriented in order to define $\star$). In this setting, the metric $g$ and the 2-form $F$ have to satisfy the coupled system
\begin{align*}
	Ric_0 &= - [F\circ F]_0,\\
	dF &=0,\\
	d\star F &=0,
\end{align*}
known as the \emph{Einstein--Maxwell Equations}. Here $[F\circ F]_0 = F_{is}{F^{s}}_j - \frac{1}{4} F_{st} F^{st} g_{ij}$ is the trace-free part of the composition of $F$ with itself where $F$ is thought as an endomorphism of the tangent bundle after raising an index. This term (up to a constant) is what physicists call the \emph{stress-energy-tensor} of the electro-magnetic field. 

Although Einstein--Maxwell Equations can be considered in any dimension $n\geq 4$, the $4$-dimensional case has a privileged status, because in this dimension, the  equations imply that the solutions must have constant-scalar-curvature~\cite{MR0125546,MR2681684}. Also, in dimension 4, if $(g,F)$ is a solution of Einstein--Maxwell Equations and $g$ is not Einstein, then $F$ is determined uniquely up to a constant: $F:=cF^+ + \frac{1}{c} F^-$, where $F^\pm=\frac{1}{2}\left(F \pm \star F\right)$ are the self-dual and the anti-self-dual parts of $F$~\cite{MR3500251}. Therefore, the Einstein--Maxwell equations can actually be thought of as having one unknown: the metric. We say that a metric is \emph{Einstein-Maxwell metric} if there is a 2-form $F$ so that $(g,F)$ is a solution of Einstein--Maxwell equations.

The Einstein--Maxwell equations also have some remarkable ties to K\"ahler Geometry. First, any K\"ahler metric with constant-scalar-curvature (cscK for short) is an Einstein--Maxwell metric. Indeed, as LeBrun~\cite{MR2681684} observed, for a cscK metric, the 2-form $F$ can be chosen as $F =   \frac{1}{2} \omega +\rho_0 $ where $\omega$ is the K\"ahler form, and $\rho_0 := Ric_0(J\cdot,\cdot)$ is the trace-free Ricci form of the metric. Second, more generally, a constant-scalar-curvature metric in a conformal class of a K\"ahler metric is Einstein--Maxwell if the conformal factor is a \emph{holomorphy potential}~\cite{MR3574879,MR3327057}. These observations lead to many examples of Einstein--Maxwell Metrics. Any cscK metric on a complex surface is a solution. Recently, some conformally K\"ahler solutions have been discovered on Hirzebruch Surfaces and more generally on so-called minimal Ruled Surfaces fibered over Riemann Surfaces of any genus~\cite{MR2681684,MR3327057,MR3500251,MR3521556,Apostolov:2015aa}. We also
refer the reader to~\cite{Shao:2013aa,MR3897025,Lahdili:2017ab,MR3868215,Futaki:2017ab,Futaki:2018aa} for more about obstructions to the existence of Einstein--Maxwell metrics.

In this paper, in pursuit of finding new examples, we look into $4$-dimensional Lie algebras. The $4$-dimensional Lie algebras were already classified by Mubarakzyanov~\cite{MR0155871} (a list can be found in~\cite{alg_class}), and the (automorphism-reduced) form of left-invariant metrics on these algebras were computed by Karki in his thesis~\cite{MR3507749}, where he also determined all left-invariant \emph{Einstein metrics} on $4$-dimensional Lie algebras up to automorphisms of the Lie algebra. Here, we find the left-invariant \emph{non-Einstein} solutions to the Einstein--Maxwell equations (up to automorphisms):

\begin{thm}\label{main}
The $4$-dimensional Lie algebras admitting left-invariant non-Einstein solutions to the Einstein--Maxwell equations are:
\begin{enumerate}[label=\arabic*)]
\item $2\mathcal{A}_2:\quad\,\,\,\,\,\,\,\,\,\,[e_1,e_2]= e_2,\quad[e_3,e_4]= e_4.$ \\
\item $\mathcal{A}_2\oplus 2\mathcal{A}_1:\,[e_1,e_2]= e_2.$\\
\item $\mathcal{A}_{4,6}^{a,0}:\quad\quad\,\, [e_1,e_4]= ae_1,\quad [e_2,e_4]=-e_3,\quad [e_3,e_4]=e_2,$ with $a\neq 0.$\\
\item $\mathcal{A}_{4,9}^{-\frac{1}{2}}:\quad\quad \,[e_2,e_3]=e_1,\quad [e_1,e_4]=\frac{1}{2}e_1,\quad [e_2,e_4]=e_2,\quad [e_3,e_4]=-\frac{1}{2}e_3.$
\end{enumerate}

\end{thm}

Here, we use the same notation of Lie algebras as in ~\cite{alg_class}. These solutions turn out to be K\"ahler with the fixed orientation $e^1\wedge e^2\wedge e^3\wedge e^4$
except on $\mathcal{A}_{4,9}^{-\frac{1}{2}}$ which admits a solution metric that cannot be (left-invariant) K\"ahler with the fixed orientation (however it is K\"ahler for the reverse orientation). That solution is actually a non-K\"ahler almost-K\"ahler metric (so the almost-complex structure $J$ is non-integrable) with $J$-invariant Ricci tensor (described in~\cite{MR2148908}). Indeed a non-K\"ahler almost-K\"ahler metric metric with $J$-invariant Ricci tensor of constant scalar curvature is a solution to the Einstein--Maxwell equations (because the Ricci form is closed in that case~\cite{MR1317012} hence the same argument for cscK metrics applies).

We also remark that $2\mathcal{A}_2$ is the only algebra
which admits a left-invariant Eintein metric and also a non-Einstein solution to the Einstein--Maxwell equations.
Furthermore, we remark that the corresponding Lie groups to all these Lie algebras admit no compact quotient.

\subsection*{Acknowledgements} The authors are very thankful to Luigi Vezzoni for his valuable help. Both authors are supported in part by a PSC-CUNY research award \#61768-00 49. 

%\section{Preliminaries}

\section{Left-invariant non-Einstein solutions to the Einstein--Maxwell equations}
We present in this section the list of all $4$-dimensional Lie algebras admitting non-Einstein solutions to the Maxwell--Einstein equations. We
give an explicit description of the solutions up to the automorphisms of the Lie algebra. In order to do so, we went over the list of  $4$-dimensional Lie algebras in~\cite{alg_class} and their (automorphism-reduced) left-invariant Riemannian metrics (as in~\cite{MR3507749})
and then used a \textit{Maple} program to determine solutions of the Einstein--Maxwell equations.
\subsection{The Lie algebra $2\mathcal{A}_2$}
The structure equations of the Lie algebra $2\mathcal{A}_2$ is
\begin{equation*}
[e_1,e_2]= e_2,\quad[e_3,e_4]= e_4.
\end{equation*}
where $\{e_i\}$ is a basis of $2\mathcal{A}_2$. This Lie algebra is not unimodular, so it does not admit a compact quotient.
Up to automorphisms of the Lie algebra (and scaling), a left-invariant metric $g$ is given by
$$g=\left[\begin{array}{cccc}1 & 0 & a_1 & a_2 \\0 & 1 & a_3 & a_4 \\a_1 & a_3 & a_5& 0 \\a_2 & a_4&0 & 1\end{array}\right],$$
where $a_i$ are constants satisfying the conditions $a_5-a_3^2-a_1^2>0$ and $(a_4^2-1)a_1^2-2a_1a_2a_3a_4+(a_3^2-a_5)a_2^2-a_4^2a_5-a_3^2+a_5>0.$
Suppose that $F=\displaystyle\sum_{1\leq i<j\leq 4}a_{ij}e^{ij}$ is a $2$-form where $a_{ij}$ are constants and $\{e^i\}$ is the dual basis of $\{e_i\}$ and $e^{ij}=e^i\wedge e^j$. Noting that $de^i = -e^i[e_j,e_k]$,
the condition $dF=0$ implies that $a_{14}=a_{23}=a_{24}=0.$ Suppose that the orientation is $e^{1234}=e^1\wedge e^2\wedge e^3\wedge e^4$. 
Then, we have
\begin{eqnarray*}
\star F&=&\frac{1}{\sqrt{\det g}}\left(a_4a_1a_{12}-a_2a_3a_{12}+a_2a_{13}+a_{34}\right)e^{12}-\frac{1}{\sqrt{\det g}}\left(a_2a_5a_{12}--a_2a_3a_{13}-a_3a_{34}\right)e^{13}\\
&+&\frac{1}{\sqrt{\det g}}\left(a_2a_4a_{13}+a_1a_{12}+a_4a_{34}\right)e^{14}-\frac{1}{\sqrt{\det g}}\left(a_4a_5a_{12}-a_4a_3a_{13}+a_1a_{34}\right)e^{23}\\
&+&\frac{1}{\sqrt{\det g}}\left(a_4^2a_{13}+a_3a_{12}-a_2a_{34}-a_{13}\right)e^{24}+\frac{1}{\sqrt{\det g}}\left(a_4a_1a_{34}-a_3a_2a_{34}+a_5a_{12}-a_3a_{13}\right)e^{34}.
\end{eqnarray*}
Then, $d\star F=0$ implies the system of equations
\begin{align*} a_2a_4a_{13}+a_1a_{12}+a_4a_{34}&=0\\ a_4a_5a_{12}-a_4a_3a_{13}+a_1a_{34}&=0\\ a_4^2a_{13}+a_3a_{12}-a_2a_{34}-a_{13}&=0
\end{align*}
which have the following \textit{non-trivial} solutions (i.e. $F\neq 0$):
\begin{enumerate}
\item $a_1=a_4=0,a_3=\frac{a_2a_{34}+a_{13}}{a_{12}}$, with $a_{12}\neq 0.$ Then $[F\circ F]_0(e_3,e_4)=0$. On the other hand, the trace-free Ricci tensor $Ric_0(e_3,e_4)=-\frac{a_2a_5a_{12}^2}{2\det g}.$ Thus $a_2=0.$ Then,
$[F\circ F]_0(e_1,e_3)=0$ while $Ric_0(e_1,e_3)=-\frac{a_{13}^2}{2\det g}$ hence $a_{13}=0.$ 
We obtain then the following solution to the Einstein--Maxwell equations given by the metric
\begin{equation}\label{metric1}
g=e^1\otimes e^1+e^2\otimes e^2+a_5\,e^3\otimes e^3+e^4\otimes e^4,
\end{equation}
and $F=a_{12}e^{12}+a_{34}e^{34}$ such that 
\begin{equation}\label{relation1}
a_5=\frac{1+a_{34}^2}{1+a_{12}^2}\neq 1.
\end{equation}
Actually, when $a_5=1$ in~(\ref{metric1}), the metric $g$ is then Einstein. In fact, the trace-free Ricci tensor of $g$ is given by
\begin{equation*}
Ric_0=\left(\frac{1-a_5}{2a_5}\right)e^1\otimes e^1+\left(\frac{1-a_5}{2a_5}\right)e^2\otimes e^2+\left(\frac{a_5-1}{2}\right)e^3\otimes e^3+\left(\frac{a_5-1}{2a_5}\right)e^4\otimes e^4.
\end{equation*}
Moreover, we have
\[
F^\pm=\frac{1}{2}\left(\pm\frac{a_{34}}{\sqrt{a_5}}+a_{12}\right)e^{12}+\frac{1}{2}\left(\pm\sqrt{a_5}a_{12}+a_{34}\right)e^{34}.
\]
%and 
%\[
%F^-=\frac{1}{2}\left(-\frac{a_{34}}{\sqrt{a_5}}+a_{12}\right)e^{12}+\frac{1}{2}\left(-\sqrt{a_5}a_{12}+a_{34}\right)e^{34}.
%\]
Furthermore, the metric $g$ is K\"ahler with respect to the K\"ahler form
\[
\omega=e^{12}+\sqrt{a_5}e^{34},
\]
with the trace-free Ricci form given by
\[
\rho_0=\frac{1-a_5}{a_5}e^{12}+\frac{a_5-1}{\sqrt{a_5}}e^{34}.
\]
Using the relation~(\ref{relation1}), we have then
\[
\frac{1}{2}\omega=\frac{1}{\left(\frac{a_{34}}{\sqrt{a_5}}+a_{12}\right)}F^+,\quad \rho_0=\left(\frac{a_{34}}{\sqrt{a_5}}+a_{12}\right)F^-.
\]
\item $a_{12}=0,a_{1}=0,a_4=0,a_{13}=-a_2a_{34}.$ Then $[F\circ F]_0(e_1,e_2)=[F\circ F]_0(e_3,e_4)=0$ implies that $a_2=a_3=0$ and we get again the solution~(\ref{metric1})
(with $a_{12}=0$).
\item $a_1=-\frac{a_4(_{13}a_2+a_{34})}{a_{12}},a_3=\frac{a_{13}+a_2a_{34}-a_{13}a_4^2}{a_{12}},a_5=\frac{a_{34}^2+a_{13}^2+2a_{13}a_{34}a_2-a_{13}^2a_4^2}{a_{12}^2}.$
Then $[F\circ F]_0\equiv 0$; so any Einstein--Maxwell metric has to be Einstein.
\end{enumerate}

%Then $d\ast F=0$ implies that
%\begin{enumerate}
%\item $a_{12}=0,a_{1}=0,a_4=0,a_{13}=-a_2a_{34}.$ Then $K_{12}=K_{34}=0$ implies that $a_2=a_3=0$ and we get a solution given by $a_5=a_{34}^2+1$ so $F=a_{34}e^{34}$ ($g$ is K\"ahler with respect to $\omega=e^{12}+\tilde{a}_{34}e^{34}$, where $\tilde{a}_{34}^2=a_5$). If $a_{34}=0$ then we get a K\"ahler-Einstein solution.\\
%%\todo{arrived here}
%\item $a_1=-\frac{a_4(_{13}a_2+a_{34})}{a_{12}},a_3=\frac{a_{13}+a_2a_{34}-a_{13}a_4^2}{a_{12}},a_5=\frac{a_{34}^2+a_{13}^2+2a_{13}a_{34}a_2-a_{13}^2a_4^2}{a_{12}^2}$
%Then $K\equiv 0.$\\
%\item $a_1=0,a_4=0,a_3=\frac{a_2a_{34}+a_{13}}{a_{12}}.$ Then $K_{34}=K_{13}=0$ implies $a_{13}=a_2=0.$ Then a solution to Einstein-Maxwell
%equation is given by $a_5=\frac{1+a_{34}^2}{1+a_{12}^2}$ ($a_{12}\neq 0$) so $F=a_{12}e^{12}+a_{34}e^{34}.$  
%\end{enumerate}
%
\subsection{The Lie algebra $\mathcal{A}_2\oplus 2\mathcal{A}_1$}
The structure equation is $[e_1,e_2]= e_2$. This Lie algebra is not unimodular, so it does not admit a compact quotient. Moreover, it does not admit any left-invariant Einstein metric~\cite{MR3507749}.
Up to automorphisms of the Lie algebra (and scaling), a left-invariant metric $g$ is given by
$$g=\left[\begin{array}{cccc}1 & 0 & 0 & 0 \\0 & 1 & a_1 & a_2 \\0 & a_1 & 1& 0 \\0 & a_2&0 & 1\end{array}\right],$$
with the conditions $1-a_1^2>0$ and $1-a_1^2-a_2^2>0.$
%Suppose that $F=\displaystyle\sum_{1\leq i<j\leq 4}a_{ij}e^{ij}.$ Then,
The condition $dF=0$ implies that $a_{23}=a_{24}=0.$ 
%Suppose that the orientation is $e^{1234}$. 
The condition $d\star F=0$ implies the following non-trivial solutions:
\begin{enumerate}
\item $a_{13}=a_1a_{12},a_{14}=a_2=0$. To get a solution to the Einstein--Maxwell equations we need $a_1=0$ and so a solution is given by
\begin{equation*}
g=e^1\otimes e^1+e^2\otimes e^2+\,e^3\otimes e^3+e^4\otimes e^4,
\end{equation*}
and $F=a_{12}e^{12}+a_{34}e^{34}$ such that 
\begin{equation}\label{relation2}
a_{34}^2-a_{12}^2=1.
\end{equation}
Moreover, we have
\[
F^\pm=\frac{1}{2}\left(a_{12}\pm a_{34}\right)e^{12}+\frac{1}{2}\left(\pm a_{12}+a_{34}\right)e^{34}.
\]
%and 
%\[
%F^-=\frac{1}{2}\left(a_{12}-a_{34}\right)e^{12}+\frac{1}{2}\left(-a_{12}+a_{34}\right)e^{34}.
%\]
Furthermore, the metric $g$ is K\"ahler with respect to the K\"ahler form
\[
\omega=e^{12}+e^{34},
\]
with the trace-free Ricci form given by
\[
\rho_0=-\frac{1}{2}e^{12}+\frac{1}{2}e^{34}.
\]
Using the relation~(\ref{relation2}), we have then
\[
\frac{1}{2}\omega=\frac{1}{\left(a_{12}+a_{34}\right)}F^+,\quad \rho_0=\left(a_{12}+a_{34}\right)F^-.\\
\]
\item  $a_{12}=\frac{a_{14}}{a_2},a_{13}=\frac{a_1a_{14}}{a_2}$, with $a_2\neq 0$. Then \textit{Maple} shows that there is no solution to the Maxwell--Einstein equations.
\end{enumerate}

\subsection{The Lie algebra $\mathcal{A}_{4,6}^{a,0}$}
The structure equations are
\begin{equation*}
[e_1,e_4]= ae_1,\quad [e_2,e_4]=-e_3,\quad [e_3,e_4]=e_2,
\end{equation*}
with $a\neq 0$. This Lie algebra does not admit a compact quotient and does not admit any left-invariant Einstein metric~\cite{alg_class}.
Up to automorphisms of the Lie algebra, a left-invariant metric $g$ is given by
$$g=\left[\begin{array}{cccc}1 & a_1 & a_2 & 0 \\a_1 & 1 & 0 & 0 \\a_2 & 0 & a_3& 0 \\0 & 0 & 0 & 1\end{array}\right],$$
with the conditions $a_3-a_3a_1^2-a_2^2>0$ and $1-a_1^2>0.$
%Suppose that $F=\displaystyle\sum_{1\leq i<j\leq 4}a_{ij}e^{ij}.$ 
The condition $dF=0$ implies that $a_{12}=a_{13}=0.$
%Suppose that the orientation is $e^{1234}$. 
Then the condition $d\star F=0$ implies
\[a_1^2a_{34}-a_2a_1a_{24}+a_2a_{14}-a_{34}=0,  \]
\[a_3a_1a_{14}-a_1a_2a_{34}+a_2^2a_{24}-a_3a_{24}=0.\]
Then, we distinguish two cases:
\begin{enumerate}
\item $a_2=0,a_{34}=0,a_{24}=a_1a_{14}$. To get a solution to the Einstein--Maxwell equations, we need
\[
a_1=0,\qquad a_3=1,\qquad a_{23}^2-a_{14}^2=a^2.
\]
Then, the Einstein--Maxwell metric is
\begin{equation*}
g=e^1\otimes e^1+e^2\otimes e^2+\,e^3\otimes e^3+e^4\otimes e^4,
\end{equation*}
and $F=a_{14}e^{14}+a_{23}e^{23}$ such that
\begin{equation}\label{relation3}
a_{23}^2-a_{14}^2=a^2,
\end{equation}
with $a\neq 0.$ If $a=0$, then $g$ is Einstein. In addition, we have
\begin{align*}
F^+&=\frac{1}{2}\left(a_{14}+a_{23}\right)e^{14}+\frac{1}{2}\left(a_{14}+a_{23}\right)e^{23}\\
F^-&=\frac{1}{2}\left(a_{14}-a_{23}\right)e^{12}+\frac{1}{2}\left(-a_{14}+a_{23}\right)e^{34}.
\end{align*}
Furthermore, the metric $g$ is K\"ahler with respect to the K\"ahler form
\[
\omega=e^{14}+e^{23}
\]
with the trace-free Ricci form given by
\[
\rho_0=-\frac{a^2}{2}e^{14}+\frac{a^2}{2}e^{34}.
\]
Using the relation~(\ref{relation3}), we have then
\[
\frac{1}{2}\omega=\frac{1}{\left(a_{14}+a_{23}\right)}F^+,\quad \rho_0=\left(a_{14}+a_{23}\right)F^-.
\]

\item $a_2\neq 0$, $a_{14}=\frac{a_{34}}{a_2}$ and $a_{24}=\frac{a_1a_{34}}{a_2}.$ \textit{Maple} shows that there is no non-Einstein solution to the Einstein--Maxwell equations.

\end{enumerate}
\subsection{The Lie algebra $\mathcal{A}_{4,9}^{-\frac{1}{2}}$}
The structure equations of the Lie algebra $\mathcal{A}_{4,9}^{-\frac{1}{2}}$ are
\begin{equation*}
[e_2,e_3]=e_1,\quad [e_1,e_4]=\frac{1}{2}e_1,\quad [e_2,e_4]=e_2,\quad [e_3,e_4]=-\frac{1}{2}e_3.
\end{equation*}
This Lie algebra does not admit a compact quotient and does not admit any left-invariant Einstein metric~\cite{alg_class}.
Up to automorphisms of the Lie algebra, a left-invariant metric $g$ is given by
$$g=\left[\begin{array}{cccc}a_1 & 0 & 0 & 0 \\0 & 1 & a_2 & a_3 \\0 & a_2 & 1& a_4 \\0 & a_3 & a_4 & 1\end{array}\right],$$
with the conditions $a_1>0,$ $1-a_2^2>0$, $1-a_4^2+2a_2a_3a_4-a_2^2-a_3^2>0.$
%Suppose that $F=\displaystyle\sum_{1\leq i<j\leq 4}a_{ij}e^{ij}.$ 
The condition $dF=0$ implies that $a_{12}=0$ and $a_{14}=\frac{1}{2}a_{23}.$
%Suppose that the orientation is $e^{1234}$. 
Then we have two cases: 
\begin{enumerate}
\item If we suppose that $a_2a_3\neq a_4$ then the condition $d\star F=0$ implies that
\[
a_{13}=\frac{4a_1a_2a_{24}a_3-4a_1a_{23}a_3^2-4a_1a_{24}a_4+a_2^2a_{23}+4a_1a_{23}-a_{23}}{a_2a_3-a_4},
\]
\[
a_{34}=a_2a_{24}-a_{23}a_3.
\]
Then \textit{Maple} shows that there is no solution to the Einstein--Maxwell equations.
\item If we suppose that $a_4=a_2a_3$ then there are two solutions to $d\star F=0$:
\subsubsection{First solution:} $a_{23}=0$ and $a_{34}=a_2a_{24}$. To get a solution to the Einstein--Maxwell equations, we need $a_1=1$ and $a_2=a_3=0$
Then, the Einstein--Maxwell metric is
\begin{equation*}
g=e^1\otimes e^1+e^2\otimes e^2+\,e^3\otimes e^3+e^4\otimes e^4
\end{equation*}
with $F=a_{13}e^{13}+a_{24}e^{24}$ such that
\begin{equation}\label{relation4}
a_{13}^2-a_{24}^2=\frac{3}{2}.
\end{equation}
It turns out that the metric $g$ is non-K\"ahler almost-K\"ahler with the orientation $e^{1234}.$ Indeed, $g$ is compatible with the closed $2$-form $
\omega=e^{13}-e^{24}$ 
inducing a non-integrable almost-complex structure $J$ defined by $Je_1=e_3$ and $Je_2=-e_4$. Moreover, its Ricci tensor is $J$-invariant. Indeed, its trace-free Ricci form is given by
$\rho_0=\frac{3}{4}e^{13}+\frac{3}{4}e^{24}$.
We have then
\begin{align*}
F^+&=\frac{1}{2}\left(a_{13}-a_{24}\right)e^{13}+\frac{1}{2}\left(a_{24}-a_{13}\right)e^{24}\\
F^-&=\frac{1}{2}\left(a_{13}+a_{24}\right)e^{13}+\frac{1}{2}\left(a_{24}+a_{13}\right)e^{24}.
\end{align*}
Furthermore, Using the relation~(\ref{relation4}), we have
\[
\frac{1}{2}\omega=\frac{2}{3}\left(a_{13}+a_{24}\right)F^+,\quad \rho_0=\frac{3}{2\left(a_{13}+a_{24}\right)}F^-.
\]
If we reverse the orientation to be $-e^{1234}$, then
\begin{align*}
F^+&=\frac{1}{2}\left(a_{13}+a_{24}\right)e^{13}+\frac{1}{2}\left(a_{24}+a_{13}\right)e^{24}\\
F^-&=\frac{1}{2}\left(a_{13}-a_{24}\right)e^{13}+\frac{1}{2}\left(a_{24}-a_{13}\right)e^{24}.
\end{align*}
Furthermore, the metric $g$ is K\"ahler with respect to the K\"ahler form $\omega=e_{13}+e_{24}$ 
with the trace-free Ricci form given by
$
\rho_0=\frac{3}{4}e^{13}-\frac{3}{4}e^{24}$.
So using the relation~(\ref{relation4}), we have then
\[
\frac{1}{2}\omega=\frac{1}{\left(a_{13}+a_{24}\right)}F^+,\quad \rho_0=\left(a_{13}+a_{24}\right)F^-.
\]
\subsubsection{Second solution:} $a_{34}=a_2a_{24}-a_3a_{23},a_1=\frac{1-a_2^2}{4(1-a_3^2)}$ (with $a_3\neq \pm 1$). Then \textit{Maple} shows that there is no solution to the Einstein--Maxwell equations.
%\end{enumerate}
\end{enumerate}
\section{Non-existence of Einstein--Maxwell metrics}
In this section, we will explain briefly why all the other Lie algebras do not admit any non-Einstein Einstein--Maxwell metrics.
\subsection{The Lie algebra $\mathcal{A}_{4,1}$}
The structure of the Lie algebra is 
\begin{equation*}
[e_2,e_4]= e_1,\quad [e_3,e_4]=e_2.
\end{equation*} 
Up to automorphisms of the Lie algebra, a metric $g$ is given by
$$g=\left[\begin{array}{cccc}1 & 0 & 0 & 0 \\0 & 1 & a_1 & 0 \\0 & a_1 & a_2 & 0 \\0 & 0 & 0 & 1\end{array}\right],$$
with the condition that $a_2-a_1^2>0.$
A form $F$ satisfies $dF=0$ if $$F=a_{14}e^{14}+a_{23}e^{23}+a_{24}e^{24}+a_{34}e^{34}.$$
Then,
%$$\star F=\frac{1}{\sqrt{\det g}} \left(\left(a_{34}-a_1a_{24}\right)e^{12}+ \left(a_1a_{34}-a_2a_{24}\right)e^{13}+a_{23}e^{14}+ \left(a_2-a_1^2\right)a_{14}e^{23}       \right).$$
The condition $d\star F=0$ implies that $a_{34}=a_1a_{24},$ and $a_1a_{34}=a_2a_{24}.$
Since $a_2\neq 0$, we get $a_{34}\left(1-\frac{a_1^2}{a_2}\right)=0$ hence $a_{34}=a_{24}=0.$
We deduce that a solution to $dF=d\star F=0$ is given by $F=a_{14}e^{14}+a_{23}e^{23}.$
Now, the tensor $[F\circ F]_0$ satisfies $[F\circ F]_0(e_1,e_2)=0$ while the trace free part of the Ricci tensor satisfies $Ric_0(e_1,e_2)=-\frac{1}{2}\frac{a_1}{{a_2-a_1^2}}$.
Hence $a_1=0$ and so there is no solution to the Einstein--Maxwell equations.

\subsection{The Lie algebra $\mathcal{A}_{4,2}^p$}
The structure of the Lie algebra is given by
\begin{equation*}
[e_1,e_4]= pe_1,\quad [e_2,e_4]=e_2, \quad [e_3,e_4]=e_2+e_3,
\end{equation*} 
with $p\neq 0.$
Up to automorphisms of the Lie algebra, a metric $g$ is given by
$$g=\left[\begin{array}{cccc}1 & a_1 & a_2 & 0 \\a_1 & 1 & 0 & 0 \\a_2 & 0 & a_3 & 0 \\0 & 0 & 0 & 1\end{array}\right],$$
with the conditions $a_3-a_3a_1^2-a_2^2>0$ and $a_3>0.$
The equation $dF=0$ implies 
\[
a_{23}=0, \qquad
a_{12}(p+1)=0, \qquad
a_{12}+a_{13}(p+1)=0.\]
We suppose first that $p\neq -1$. Then we get $a_{12}=a_{13}=0.$ 
%Hence,
%\begin{eqnarray*}
%\star F&=&\frac{1}{\sqrt{\det g}} \left(\left(a_{34}-a_2a_{14}+a_2a_1a_{24}-a_1^2a_{34}\right)e^{12}+  \left(a_3a_1a_{14}-a_1a_2a_{34}+a_2^2a_{24}-a_3a_{24}\right)e^{13}\right.\\
%&+&\left. \left(-a_3a_1a_{24}+a_3a_{14}-a_2a_{34}\right)e^{23}   \right).
%\end{eqnarray*}
The condition $d\star F=0$ implies that 
\[
a_{34}-a_2a_{14}+a_2a_1a_{24}-a_1^2a_{34}=0, \]
\[
a_3a_1a_{14}-a_1a_2a_{34}+a_2^2a_{24}-a_3a_{24}=0, \]
\[
-a_3a_1a_{24}+a_3a_{14}-a_2a_{34}=0.\]
From the third equation we get (since $a_3>0$) $a_{14}=a_1a_{24}+\frac{a_2}{a_3}a_{34}$. Replacing it in the second equation we get that $a_{24}=0$ because $a_3-a_3a_1^2-a_2^2>0$.
Then it is easy to deduce that $a_{34}=a_{14}=0.$
We conclude that under the hypothesis $p\neq -1$, there is no non trivial $F$ satisfying $dF=d\star F=0.$

Now, we suppose that $p=-1.$ Then $dF=0$ implies that $a_{23}=a_{12}=0.$ 
%Then,
%\begin{eqnarray*}
%\star F&=&\frac{1}{\sqrt{\det g}}  \left(\left(a_{34}-a_2a_{14}+a_2a_1a_{24}-a_1^2a_{34}\right)e^{12}+  \left(a_3a_1a_{14}-a_1a_2a_{34}+a_2^2a_{24}-a_3a_{24}\right)e^{13}\right.\\
%&+&\left. \left(-a_3a_1a_{24}+a_3a_{14}-a_2a_{34}\right)e^{23}-a_1a_{13}e^{14}-a_{13}e^{24}   \right).
%\end{eqnarray*}
From $d\star F=0,$ it follows that
\[a_{34}-a_2a_{14}+a_2a_1a_{24}-a_1^2a_{34}=0,  \qquad 
-a_3a_1a_{24}+a_3a_{14}-a_2a_{34}=0.\]
We get $a_{14}=a_1a_{24}+\frac{a_2}{a_3}a_{34}$ from the second equation. Replacing it in the first we obtain that $a_{34}=0.$
A solution $F$ of $dF=d\star F=0$ is of the form $$
F=a_{13}e^{13}+a_1a_{24}e^{14}+a_{24}e^{24},$$
and then using \textit{Maple}, it turns out that there are no solutions to Einstein--Maxwell equations.

\subsection{The Lie algebra $\mathcal{A}_{4,3}$}
The structure of the Lie algebra is
\begin{equation*}
[e_1,e_4]= e_1,\quad [e_3,e_4]=e_2.
\end{equation*}
Up to automorphisms of the Lie algebra, a metric $g$ is given by
$$g=\left[\begin{array}{cccc}1 & a_1 & a_2 & 0 \\a_1 & 1 & 0 & 0 \\a_2 & 0 & a_3 & 0 \\0 & 0 & 0 & 1\end{array}\right],$$
with the conditions that $a_3-a_3a_1^2-a_2^2>0$ and $a_3>0.$
The equation $dF=0$ implies that $a_{12}=a_{13}=0.$ 
%Now,
%\begin{eqnarray*}
%\star F&=&\frac{1}{\sqrt{\det g}}    \left[(-a_1^2a_{34}+a_2a_1a_{24}-a_2a_{14}+a_{34})e^{12}+  (a_3a_1a_{14}-a_1a_2a_{34}+a_2^2a_{24}-a_3a_{24})e^{13}\right.\\
%&+&\left. (-a_3a_1a_{24}+a_3a_{14}-a_2a_{34})e^{23}+a_{23}e^{14}+a_1a_{23}e^{24}+a_2a_{23}e^{34}   \right].
%\end{eqnarray*}
From $d\star F=0,$ it follows that
\[
-a_1^2a_{34}+a_2a_1a_{24}-a_2a_{14}+a_{34}=0 , \quad 
a_3a_1a_{14}-a_1a_2a_{34}+a_2^2a_{24}-a_3a_{24}=0.\]
We deduce then that $a_{34}=a_2a_{14},a_{24}=a_1a_{14}$ and then Maple shows that there are no solutions to the Einstein--Maxwell equations.

\subsection{The Lie algebra $\mathcal{A}_{4,4}$}
The structure of the Lie algebra is
\begin{equation*}
[e_1,e_4]= e_1,\quad [e_2,e_4]=e_1+e_2,\quad [e_3,e_4]=e_2+e_3.
\end{equation*}
Up to automorphisms of the Lie algebra, a metric $g$ is given by
$$g=\left[\begin{array}{cccc}1 & 0 & 0 & 0 \\0 & a_1 & a_2 & 0 \\0 & a_2 & a_3 & 0 \\0 & 0 & 0 & 1\end{array}\right],$$
with the conditions that $a_1a_3-a_2^2>0$ and $a_1>0.$
Now, $dF=0$ implies that $a_{12}=a_{13}=a_{23}=0.$ 
%On the other hand,
%\begin{eqnarray*}
%\star F&=&\frac{1}{\sqrt{\det g}}  \left[(a_1a_{34}-a_2a_{24})e^{12}+  (a_2a_{34}-a_3a_{24})e^{13}\right.\\
%&+&\left. a_{14}(a_1a_3-a_2^2)e^{23}   \right].
%\end{eqnarray*}
From $d\star F=0,$ it follows that
\[
 a_1a_{34}-a_2a_{24}=0,\quad 
 a_2a_{34}-a_3a_{24}=0,\quad
 a_{14}(a_1a_3-a_2^2)=0.\]
Hence $a_{14}=a_{24}=a_{34}=0$ and thus there is no non trivial solution $F.$

\subsection{The Lie algebra $\mathcal{A}_{4,5}^{a,b}$}
The structure of the Lie algebra is
\begin{equation*}
[e_1,e_4]= e_1,\quad [e_2,e_4]=ae_2,\quad [e_3,e_4]=be_3,
\end{equation*}
with $ab\neq 0,\,-1\leq a\leq b\leq 1.$
Up to automorphisms of the Lie algebra, a metric $g$ is given by
$$g=\left[\begin{array}{cccc}1 & a_1 & a_2 & 0 \\a_1 & 1 & a_3 & 0 \\a_2 & a_3 & 1 & 0 \\0 & 0 & 0 & 1\end{array}\right],$$
with the conditions that $1-a_1^2-a_2^2-a_3^2+2a_1a_2a_3>0,$ $1-a_1^2>0,$ $(1-a_2^2)(1-a_3^2)>0.$
Now, $dF=0$ implies the following solutions depending on $a$ and $b$:
\subsubsection{Case $a\neq -1,$ $b\neq -1$ and $a\neq -b$}
In this case, we have $a_{12}=a_{13}=a_{23}=0.$ The condition $d\star F=0$ implies that
\[
 a_1^2a_{34}-a_3a_1a_{14}-a_1a_2a_{24}+a_2a_{14}+a_3a_{34}-a_{34}=0,  \]
 \[
  -a_1a_2a_{34}-a_2a_3a_{14}+a_2^2a_{24}+a_1a_{14}+a_3a_{34}-a_{24}=0,\]
  \[
  -a_3a_1a_{34}+a_3^2a_{14}-a_3a_2a_{24}+a_1a_{24}+a_2a_{34}-a_{14}=0.\]
It turns out that there is no non trivial solution $F$.
%
%$$\left\{\begin{array}{c}  a_1=a_2a_3+\sqrt{1+a_2^2a_3^3-a_2^2-a_3^2},\\a_{14}=\frac{1}{1-a_3^2}\left(-a_3a_{34}\left(a_2a_3+\sqrt{1+a_2^2a_3^3-a_2^2-a_3^2}\right)-a_2a_3a_{24}+a_{24}\left(  a_2a_3+\sqrt{1+a_2^2a_3^3-a_2^2-a_3^2}\right)+a_2a_{34}   \right) \end{array}\right.$$
%and
%$$\left\{\begin{array}{c}  a_1=a_2a_3-\sqrt{1+a_2^2a_3^3-a_2^2-a_3^2},\\a_{14}=\frac{1}{1-a_3^2}\left(-a_3a_{34}\left(a_2a_3-\sqrt{1+a_2^2a_3^3-a_2^2-a_3^2}\right)-a_2a_3a_{24}+a_{24}\left(  a_2a_3-\sqrt{1+a_2^2a_3^3-a_2^2-a_3^2}\right)+a_2a_{34}   \right) \end{array}\right.$$
\subsubsection{Case $a\neq -1$ $b\neq -1$ and $a= -b$}
In this case, we have $a_{12}=a_{13}=0.$ The condition $d\star F=0$ implies that
\[
a_1^2a_{34}-a_3a_1a_{14}-a_1a_2a_{24}+a_2a_{14}+a_3a_{34}-a_{34}=0,  \]
\[
  -a_1a_2a_{34}-a_2a_3a_{14}+a_2^2a_{24}+a_1a_{14}+a_3a_{34}-a_{24}=0.\]
Then the non trivial solution is $a_{24}=a_1a_{14}$ and $a_{34}=a_2a_{14}.$
Hence, $$F=a_{14}e^{14}+a_{23}e^{23}+a_1a_{14}e^{24}+a_2a_{14}e^{34}$$
%Then we have that $K:=(F\circ F)_0$ satisfies $K_{11}=K_{44}$ then $(r_0)_{11}=(r_0)_{44}$ implies that $b=0$ or $b=\frac{a(a_2^2-a_1^2)}{a_1^2+a_2^2+4}.$
Then \textit{Maple} shows that there is no solution to the Einstein--Maxwell equations.
\subsubsection{Case $a\neq -1,$ $b=-1$ and $a\neq -b$} In this case, $a_{12}=a_{23}=0$. Then $d\star F=0$ implies $a_{14}=a_1a_{24},a_{34}=a_3a_{24}$
and it turns out that there is no solution of the Einstein--Maxwell equations.

\subsubsection{Case $a= -1,$ $b\neq -1$ and $a\neq -b$} In this case, $a_{13}=a_{23}=0$. Then $d\star F=0$ implies $a_{14}=a_2a_{34},a_{24}=a_3a_{34}$
and it turns out that there is no solution of the Einstein--Maxwell equations.
\subsubsection{Case $a= -1,$ $b= -1$} So $dF=0$ implies $a_{23}=0$. The condition $d\star F=0$ implies
$a_{14}=(-a_3a_1a_{34}-a_3a_2a_{24}+a_1a_{24}+a_2a_{34})/(1-a_3^2)$,
and it turns out that there is no solution of the Einstein--Maxwell equations.
\subsubsection{Case $a= -1,b=1$} The condition $dF=0$ implies $a_{13}=0$. The condition $d\star F=0$ implies $
a_{24}=(-a_1a_2a_{34}-a_2a_3a_{14}+a_1a_{14}+a_3a_{34})/(1-a_2^2)
$,
and it turns out that there is no solution of the Einstein--Maxwell equations.
\subsubsection{Case $a=1,b= -1,$} Then $dF=0$ implies $a_{12}=0$. The condition $d\star F=0$ implies
$a_{34}=(-a_3a_1a_{14}-a_2a_1a_{24}+a_2a_{14}+a_3a_{24})/(1-a_1^2)
$,
and it turns out that there is no solution of the Einstein--Maxwell equations.
\subsection{The Lie algebra $\mathcal{A}_{4,6}^{a,b}$}
The structure of the Lie algebra is
\begin{equation*}
[e_1,e_4]= ae_1,\quad [e_2,e_4]=be_2-e_3,\quad [e_3,e_4]=e_2+be_3,
\end{equation*}
with $a\neq 0,\, b>0.$
Up to automorphisms of the Lie algebra, a metric $g$ is given by
$$g=\left[\begin{array}{cccc}1 & a_1 & a_2 & 0 \\a_1 & 1 & 0 & 0 \\a_2 & 0 & a_3& 0 \\0 & 0 & 0 & 1\end{array}\right],$$
with the conditions that $a_3-a_3a_1^2-a_2^2>0,$ $1-a_1^2>0.$
The condition $dF=0$ implies that
\[
a_{12}(a+b)=a_{13},\quad
a_{13}(a+b)=-a_{12},\quad
b\,a_{23}=0.
\]
This implies that $a_{12}=a_{13}=0=a_{23}=0$. The condition $d\star F=0$ implies that
$$
 a_1^2a_{34}-a_2a_1a_{24}+a_2a_{14}-a_{34}=0, \quad a_3a_1a_{24}-a_3a_{14}+a_2a_{34}=0,$$
$$
 a_3a_1a_{14}-a_1a_2a_{34}+a_2^2a_{24}-a_3a_{24}=0.$$
Then there is no non-trivial solution $F$.
\subsection{The Lie algebra $\mathcal{A}_{4,7}$}
The structure of the Lie algebra is
\begin{equation*}
[e_2,e_3]= e_1,\quad [e_1,e_4]=2e_1,\quad [e_2,e_4]=e_2,\quad [e_3,e_4]=e_2+e_3.
\end{equation*}
Up to automorphisms of the Lie algebra, a metric $g$ is given by
$$g=\left[\begin{array}{cccc}a_1 & a_2 & a_3 & 0 \\a_2 & a_4 & 0 & 0 \\a_3 & 0 & 1& 0 \\0 & 0 & 0 & 1\end{array}\right],$$
with the conditions that $a_1>0,$ $a_4>0$, $a_1a_4-a_2^2-a_3^2a_4>0.$ 
The condition $dF=0$ implies that $a_{12}=a_{13}=0,a_{23}=\frac{1}{2}a_{14}.$
Then $d\star F=0$ implies that 
\begin{align*}
 a_4a_1a_{34}-a_3a_4a_{14}-a_2^2a_{34}+a_3a_2a_{24}=0,\\ 
 a_2a_3a_{34}-a_3^2a_{24}+a_1a_{24}-a_2a_{14}=0,\\
-a_4a_3a_{34}+a_4a_{14}-a_2a_{24}-\frac{1}{4}a_1a_{14}=0.
\end{align*}
Then there are two possible solutions:
\begin{enumerate}
\item $a_2=a_{24}=0$ and $a_{34}=\frac{a_3a_{14}}{a_1}$ and $a_4=\frac{a_1^2}{4(a_1-a_3^2)}$.
Then the tensor $[F\circ F]_0\equiv 0$, and so any Maxwell--Einstein metric is Einstein.\\

\item $a_1=\frac{a_2a_{14}}{a_{24}},a_3=\frac{a_2a_{34}}{a_{24}}$ and $a_{4}=\frac{a_2(a_{14}^2+4a_{24}^2)}{4(-a_2a_{34}^2+a_{14}a_{24})}$ (with $a_{24}\neq 0$ and $-a_2a_{34}^2+a_{14}a_{24}\neq 0$ otherwise we are in the first case). 
We get again $[F\circ F]_0\equiv 0$.
\end{enumerate}

\subsection{The Lie algebra $\mathcal{A}_{4,8}$}
The structure of the Lie algebra is
\begin{equation*}
[e_2,e_3]= e_1,\quad [e_2,e_4]=e_2,\quad [e_3,e_4]=-e_3.
\end{equation*}
Up to automorphisms of the Lie algebra, a metric $g$ is given by
$$g=\left[\begin{array}{cccc}a_1 & 0 & 0 & 0 \\0 & 1 & a_2 & a_3 \\0 & a_2 & 1& a_4 \\0 & a_3 & a_4 & 1\end{array}\right],$$
with the conditions that $a_1>0,$ $1-a_2^2>0$, $1-a_4^2+2a_2a_3a_4-a_2^2-a_3^2>0,$ 
The conditions $dF=0$ implies that $a_{12}=a_{13}=a_{14}=0.$
However, there is no non trivial to the equation $d\star F=0.$

\subsection{The Lie algebra $\mathcal{A}_{4,9}^b$}
The structure of the Lie algebra is
\begin{equation*}
[e_2,e_3]= e_1,\quad [e_1,e_4]=(b+1) e_1,\quad [e_2,e_4]=e_2,\quad [e_3,e_4]=be_3,
\end{equation*}
with the conditions $-1<b\leq 1$ and $b\neq -\frac{1}{2}.$
Up to automorphisms of the Lie algebra, a metric $g$ is given by
$$g=\left[\begin{array}{cccc}a_1 & 0 & 0 & 0 \\0 & 1 & a_2 & a_3 \\0 & a_2 & 1& a_4 \\0 & a_3 & a_4 & 1\end{array}\right],$$
with the conditions $a_1>0,$ $1-a_2^2>0$, $1-a_4^2+2a_2a_3a_4-a_2^2-a_3^2>0.$
From the condition $dF=0$, we get $a_{12}=0,a_{13}=0, a_{14}=(1+b)a_{23}.$ Then $d\star F=0$ implies that ($a_{23}\neq 0$ otherwise $F$ is trivial)
\[
  a_3=\frac{a_2^2a_{34}+a_2a_{23}a_4-a_{34}}{a_{23}},\qquad 
      a_{24}=a_2a_{34}+a_{23}a_4,\]
       \[
b=\frac{a_1a_2^2a_{34}^2-a_1a_{23}^2a_4^2+a_2^2a_{23}^2+a_1a_{23}^2-a_1a_{34}^2-a_{23}^2}{a_{23}^2(1-a_2^2)}.\]
\textit{Maple} shows then that there is no solution to the Einstein--Maxwell equations.

\subsection{The Lie algebra $\mathcal{A}_{4,10}$}
The structure of the Lie algebra is
\begin{equation*}
[e_2,e_3]= e_1,\quad [e_2,e_4]=- e_3,\quad [e_3,e_4]=e_2,
\end{equation*}
Up to automorphisms of the Lie algebra, a metric $g$ is given by
$$g=\left[\begin{array}{cccc}a_1 & 0 & 0 & 0 \\0 & a_2 & 0 & a_3 \\0 & 0 & 1& a_4 \\0 & a_3 & a_4 & 1\end{array}\right],$$
with the conditions that $a_1>0,$ $a_2>0$, $a_2-a_2a_4^2-a_3^2>0.$
Now, $dF=0$ implies that $a_{12}=a_{13}=a_{14}=0.$ Then $d\star F=0$ implies that $a_2=\frac{a_{23}^2(1-a_4^2)}{a_{34}^2}$, $ a_{24}=a_{23}a_4$, $a_3=\frac{a_{23}(a_4^2-1)}{a_{34}}$, ($a_{34}\neq 0$ otherwise $F$ is trivial). But then the determinant of $g$ is 0.

\subsection{The Lie algebra $\mathcal{A}_{4,11}^a$}
The structure of the Lie algebra is
\begin{equation*}
[e_2,e_3]= e_1,\quad[e_1,e_4]=2a e_1,\quad [e_2,e_4]=ae_2- e_3,\quad [e_3,e_4]=e_2+ae_3,
\end{equation*}
with $a>0$. Up to automorphisms of the Lie algebra, a metric $g$ is given by
$$g=\left[\begin{array}{cccc}a_1 & 0 & 0 & 0 \\0 & a_2 & 0 & a_3 \\0 & 0 & 1& a_4 \\0 & a_3 & a_4 & 1\end{array}\right],$$
with the conditions that $a_1>0,$ $a_2>0$, $a_2-a_2a_4^2-a_3^2>0.$
The condition $dF=0$ implies that $a_{12}=a_{13}=0$ and $a_{14}=2aa_{23}.$ Moreover,
The condition $d\star F=0$ implies $a_1=\frac{4a_2^2a^2}{a_2-a_2a_4^2-a_3^2}$, $a_{24}=a_{23}a_4$,
$a_{34}=-\frac{a_{23}a_3}{a_2}$. Then $[F\circ F]_0\equiv 0$ and hence any Einstein--Maxwell metric is Einstein.

\subsection{The Lie algebra $\mathcal{A}_{4,12}$}
The structure of the Lie algebra is
\begin{equation*}
[e_1,e_3]= e_1,\quad[e_2,e_3]=e_2,\quad [e_1,e_4]=-e_2,\quad [e_2,e_4]=e_1.
\end{equation*}
Up to automorphisms of the Lie algebra, a metric $g$ is given by
$$g=\left[\begin{array}{cccc}1 & 0 & 0 & 0 \\0 & a_1 & a_2 & a_3 \\0 & a_2 & 1& a_4 \\0 & a_3 & a_4 & a_5\end{array}\right],$$
with the conditions that $a_1>0,$ $a_1-a_2^2>0$, $a_1a_5-a_1a_4^2-a_2^2a_5+2a_2a_3a_4-a_3^2>0.$
Then $dF=0$ implies that $a_{12}=0,a_{13}=a_{24},a_{14}=-a_{23}.$
The condition $d\star F=0$ implies the following different solutions
\begin{enumerate}
\item $a_2=a_{23}=a_{34}=a_4=0$ and $a_5=\frac{a_3^2+1}{a_1}.$ Then $[F\circ F]_0\equiv 0$ and so any Einstein--Maxwell metric is Einstein.\\
\item $a_{24}=0,a_{34}=-\frac{a_3a_{23}}{a_1},a_4=\frac{a_2a_3}{a_1},a_5=\frac{-a_1a_2^2+a_1^2+a_3^2}{a_1}.$  Then we have $[F\circ F]_0\equiv 0$.\\
\item Suppose that $a_{24}\neq 0,a_{23}\neq -a_{24}a_4,a_{24}\neq a_{23}a_4,a_{23}\neq 0$. Then
\[
 a_1=-\frac{a_{23}(a_{23}a_4-a_{24})}{a_{24}(a_{24}a_4+a_{23})}, \qquad  a_{34}=\frac{a_3a_{24}(a_4a_{24}+a_{23})}{a_{23}a_4-a_{24}}, \qquad a_2=0,\]
   \[
a_5=-\frac{a_{24}^3a_3^2a_4+a_{23}^3a_4^2+a_{23}a_{24}^2a_3^2-a_{23}a_{24}^2a_4^2-2a_{23}^2a_{24}a_4+a_{24}^3a_4+a_{23}a_{24}^2 }{a_{23}a_{24}(a_{23}a_4-a_{24})}.\]
Then $[F\circ F]_0\equiv 0$.\\
\item If we suppose that $a_{24}=0$ then we are in the case $ii).$\\
\item If we suppose that $a_{23}=0$ then we are in the case $i).$\\
\item If we suppose that $a_{24}=a_{23}a_4$. Then either \\
\begin{enumerate}
\item $a_2=a_4=0,a_{34}=-\frac{a_3a_{23}}{a_1}$ and $a_5=\frac{a_1^2+a_3^2}{a_1}$. Then $[F\circ F]_0\equiv 0.$\\
\item $a_2\neq 0$ and then
\[
a_3=\frac{(a_1^2a_4^2-a_1a_2^2+a_1^2+a_2^2)a_4}{a_2(a_1a_4^2+1)},\qquad 
a_{34}= -\frac{a_{23}a_4(-a2^2a_4^2+a_1a_4^2-a_2^2+a_1)}{a_2(a_1a_4^2+1)},\]
\[
a_5=\frac{a_1^3a_4^6-2a_1^2a_2^2a_4^4-a_2^4a_4^4+2a_1^3a_4^4+3a_1a_2^2a_4^4-a_1^2a_2^2a_4^2-2a_2^4a_4^2+a_1^3a_4^2+2a_1a_2^2a_4^2-a_2^4+a_2^2a_4^2+a_1a_2^2}{a_2^2(a_1a_4^2+1)^2}.\]
Then $[F\circ F]_0\equiv 0.$
\end{enumerate}

\item If we suppose that $a_{23}=-a_4a_{24}$ Then either\\
\begin{enumerate}
\item $a_2=0,a_{34}=0,a_5=\frac{1+a_3^2}{a_1}$. Then $[F\circ F]_0(e_1,e_3)=0$ implies that $a_3=0$ and there will be no solution to the Einstein--Maxwell equations.\\
\item  If $a_2\neq 0$, then
\[
 a_3 =\frac{ (a_1a_2^2+a_1a_4^2-a_2^2+a_1)a_4}{a_2(a_4^2+a_1)},\qquad
a_{34} =\frac{ a_{24}(a_2^2a_4^2+a_4^4+a_2^2+a_4^2)}{a_2(a_4^2+a_1)}, \]
\[
a_5 =\frac{ -a_2^4a_4^4+3a_1a_2^2a_4^4+a_1a_4^6+a_1^2a_2^2a_4^2-2a_2^4a_4^2-2a_2^2a_4^4+2a_1a_2^2a_4^2+2a_1a_4^4-a_2^4-a_2^2a_4^2+a_1a_2^2+a_1a_4^2}{a_2^2(a_4^2+a_1)^2}.\]
Then $[F\circ F]_0\equiv 0.$
\end{enumerate}
\item $a_2\neq 0$ and
\[
a_3=\frac{a_1^2a_{24}^2a_4-a_1a_2^2a_{23}a_{24}+a_1^2a_{23}a_{24}+a_1a_{23}^2a_4+a_2^2a_{23}a_{24}-a_1a_{23}a_{24}}{a_2(a_1a_{24}^2+a_{23}^2)},\]
\[a_{34}=-\frac{a_1^2a_{24}^2a_{23}a_4-a_2^2a_{23}^2a_{24}-a_2^2a_{24}^3+a_1a_{23}^2a_{24}+a_{23}^3a_4-a_{23}^2a_{24}}{a_2(a_1a_{24}^2+a_{23}^2)},\]
\begin{eqnarray*}
a_2^2(a_1a_{24}^2+a_{23}^2)a_5&=&a_1^3a_{24}^4a_4^2-2a_1^2a_2^2a_{23}a_{24}^3a_4+2a_1^3a_{23}a_{24}^3a_4-a_1^2a_2^2a_{23}^2a_{24}^2\\
&+&2a_1^2a_{23}^2a_{24}^2a_4^2-2a_1a_2^2a_{23}^3a_{24}a_4+2a_1a_2^2a_{23}a_{24}^3a_4-a_2^4a_{23}^4-2a_2^4a_{23}^2a_{24}^2-a_2^4a_{24}^4\\
&+&a_1^3a_{23}^2a_{24}^2+2a_1^2a_{23}^3a_{24}a_4-2a_1^2a_{23}a_{24}^3a_4+a_1a_2^2a_{23}^4+4a_1a_2^2a_{23}^2a_{24}^2+a_1a_2^2a_{24}^4\\
&+&a_1a_{23}^4a_4^2+2a_2^2a_{23}^3a_{24}a_4-2a_1^2a_{23}^2a_{24}^2-2a_1a_{23}^3a_{24}a_4-a_2^2a_{23}^2a_{24}^2+a_1a_{23}^2a_{24}^2.
\end{eqnarray*}
Then $[F\circ F]_0\equiv 0.$ Here $a_{23}\neq 0$ and $a_{24}\neq 0$ otherwise $F$ vanishes.
\end{enumerate}

\subsection{The Lie algebra $\mathcal{A}_{3,1}\oplus\mathcal{A}_{1}$}
The structure of the Lie algebra is $[e_3,e_4]= e_1$.
Up to automorphisms of the Lie algebra, a metric $g$ is given by
$$g=\left[\begin{array}{cccc}a_1 & 0 & 0 & 0 \\0 & 1 & 0 & 0 \\0 & 0 & 1& 0 \\0 & 0 &0 & 1\end{array}\right],$$
with the condition $a_1>0.$ The condition $dF=0$ implies $a_{12}=0$ and then $d\star F=0$ implies $a_{34}=0.$
\textit{Maple} shows then that there is no solution to the Einstein--Maxwell equation.

\subsection{The Lie algebra $\mathcal{A}_{3,2}\oplus\mathcal{A}_{1}$}
The structure of the Lie algebra is
\begin{equation*}
[e_1,e_3]= e_1,\quad[e_2,e_3]= e_1+e_2.
\end{equation*}
Up to automorphisms of the Lie algebra, a metric $g$ is given by
$$g=\left[\begin{array}{cccc}a_1 & 0 & 0 & a_2 \\0 & 1 & 0 & a_3 \\0 & 0 & 1& 0 \\a_2 & a_3 &0 & 1\end{array}\right],$$
with the conditions $a_1>0,a_1-a_2^2-a_1a_3^2>0.$ The condition $dF=0$ implies $a_{12}=a_{14}=a_{24}=0.$
Then $d\star F=0$ has no non trivial solution.

\subsection{The Lie algebra $\mathcal{A}_{3,3}\oplus\mathcal{A}_{1}$}
The structure of the Lie algebra is
\begin{equation*}
[e_1,e_3]= e_1,\quad[e_2,e_3]= e_2.
\end{equation*}
Up to automorphisms of the Lie algebra, a metric $g$ is given by
$$g=\left[\begin{array}{cccc}1 & a_1 & 0 & a_2 \\a_1 & 1 & 0 & a_3 \\0 & 0 & 1& 0 \\a_2 & a_3 &0 & 1\end{array}\right],$$
with the conditions $1-a_1^2>0,1-a_1^2-a_2^2-a_3^2+2a_1a_2a_3>0$. The condition $dF=0$ implies $a_{12}=a_{14}=a_{24}=0.$
Then the only solution to $d\star F=0$ is the trivial solution.

\subsection{The Lie algebra $\mathcal{A}_{3,4}\oplus\mathcal{A}_{1}$}
The structure of the Lie algebra is
\begin{equation*}
[e_1,e_3]= e_1,\quad[e_2,e_3]=- e_2.
\end{equation*}
Up to automorphisms of the Lie algebra, a metric $g$ is given by
$$g=\left[\begin{array}{cccc}1 & a_1 & 0 & a_2 \\a_1 & 1 & 0 & a_3 \\0 & 0 & 1& 0 \\a_2 & a_3 &0 & 1\end{array}\right],$$
with the conditions $1-a_1^2>0,1-a_1^2-a_2^2-a_3^2+2a_1a_2a_3>0$. The condition $dF=0$ implies 
$a_{14}=a_{24}=0.$ The condition $d\star F=0$ implies that $a_{13}=-a_2a_{34},a_{23}=-a_3a_{34}.$
Then \textit{Maple} shows that there is no solution to the Einstein--Maxwell equations.

\subsection{The Lie algebra $\mathcal{A}_{3,5}\oplus\mathcal{A}_{1}$}
The structure of the Lie algebra is
\begin{equation*}
[e_1,e_3]= e_1,\quad[e_2,e_3]=a e_2,
\end{equation*}
with $0<|a|<1$. Up to automorphisms of the Lie algebra, a metric $g$ is given by
$$g=\left[\begin{array}{cccc}1 & a_1 & 0 & a_2 \\a_1 & 1 & 0 & a_3 \\0 & 0 & 1& 0 \\a_2 & a_3 &0 & 1\end{array}\right],$$
with the conditions $1-a_1^2>0,1-a_1^2-a_2^2-a_3^2+2a_1a_2a_3>0$. The condition $dF=0$ implies $a_{12}=a_{14}=a_{24}=0$
Then the only solution to $d\star F=0$ is the trivial solution.

\subsection{Lie algebra $\mathcal{A}_{3,6}\oplus\mathcal{A}_{1}$}
The structure of the Lie algebra is
\begin{equation*}
[e_1,e_3]= -e_2,\quad[e_2,e_3]= e_1.
\end{equation*}
Up to automorphisms of the Lie algebra, a metric $g$ is given by
$$g=\left[\begin{array}{cccc}1 & 0 & 0 & a_1 \\0 & a_2 & 0 & a_3 \\0 & 0 & 1& 0 \\a_1 & a_3 &0 & 1\end{array}\right],$$
with the conditions $a_2>0,a_2-a_3^2-a_2a_1^2>0$. The condition $dF=0$ implies $a_{14}=a_{24}=0.$
The condition $d\star F=0$ implies the following solutions
\begin{enumerate}
\item $a_3=0,a_{13}=-a_1a_{34},a_{23}=0.$
Then  $a_1=0,a_2=1,a_{12}=\pm a_{34}$ and hence $[F\circ F]\equiv 0$ and so any Einstein--Maxwell metric is Einstein.\\

\item $a_3\neq 0,a_{13}=\frac{a_1a_{23}}{a_3},a_{34}=-\frac{a_{23}}{a_3}$. Then \textit{Maple} shows that there is no solution to the Einstein--Maxwell equations.
%We get a {\textbf{family of left invariant non K\"ahler solutions}} to the Einstein-Maxwell equation given
%by $a_1^2=\frac{1}{2}-a_3^2,a_2=1,a_{12}=\pm \frac{1}{2}\sqrt{2+\sqrt{2}},a_{23}=\pm 2^{-\frac{1}{4}}a_3$. Here the solution is well defined for any $0<|a_3|<\frac{1}{\sqrt{2}}.$
%The metric $g$ is then
%$$g=\left[\begin{array}{cccc}1 & 0 & 0 & \pm \sqrt{\frac{1}{2}-a_3^2} \\0 & 1 & 0 & a_3 \\0 & 0 & 1& 0 \\  \pm \sqrt{\frac{1}{2}-a_3^2}  & a_3 &0 & 1\end{array}\right]$$
%and 
%\[
%F=\pm \frac{1}{2}\sqrt{2+\sqrt{2}}\,e^{12}\pm \sqrt{\frac{1}{2}-a_3^2} \, 2^{-\frac{1}{4}}e^{13}\pm 2^{-\frac{1}{4}}a_3e^{23}\mp2^{-\frac{1}{4}}e^{34}.
%\]
%A closed form $\omega$ compatible with $g$ induce one of the following almost-complex structures
%$$J=\left[\begin{array}{cccc} \sqrt{2}\,a_3a_1  &(\sqrt{2}a_3^2-\sqrt{2})&0&0 \\  \sqrt{2}(2a_3^2+1)  &   -\sqrt{2}\,a_3a_1& 0 & 0 \\- a_1& -a_3 & 0& -1 \\ -\sqrt{2}\,a_3& \sqrt{2}\,a_1 &1 & 0\end{array}\right]$$
%or
%$$J=\left[\begin{array}{cccc} \sqrt{2}\,a_3a_1  &(\sqrt{2}a_3^2-\sqrt{2})&0&0 \\  \sqrt{2}(2a_3^2+1)  &   -\sqrt{2}\,a_3a_1& 0 & 0 \\ a_1& a_3 & 0& 1 \\ -\sqrt{2}\,a_3& \sqrt{2}\,a_1 &-1 & 0\end{array}\right]$$
%However, it turns out that $J$ is not integrable. 

\end{enumerate}

\subsection{The Lie algebra $\mathcal{A}_{3,7}\oplus\mathcal{A}_{1}$}
The structure of the Lie algebra is
\begin{equation*}
[e_1,e_3]= ae_1-e_2,\quad[e_2,e_3]= e_1+ae_2,
\end{equation*}
with $a>0$. Up to automorphisms of the Lie algebra, a metric $g$ is given by
$$g=\left[\begin{array}{cccc}1 & 0 & 0 & a_1 \\0 & a_2 & 0 & a_3 \\0 & 0 & 1& 0 \\a_1 & a_3 &0 & 1\end{array}\right],$$
with the conditions $a_2>0,a_2-a_3^2-a_2a_1^2>0$. The condition $dF=0$ implies that $a_{12}=a_{14}=a_{24}=0.$
Then the only solution to $d\star F=0$ is the trivial solution.

\subsection{Lie algebra $\mathcal{A}_{3,8}\oplus\mathcal{A}_{1}$}
The structure of the Lie algebra is
\begin{equation*}
[e_1,e_3]=-2e_2,\quad [e_1,e_2]= e_1,\quad [e_2,e_3]= e_3. 
\end{equation*}
Up to automorphisms of the Lie algebra, a metric $g$ is given by
$$g=\left[\begin{array}{cccc}a_1 & 0 & 0 & a_4 \\0 & a_2 & 0 & a_5 \\0 & 0 & a_3& a_6 \\a_4 & a_5&a_6 & 1\end{array}\right],$$
with the conditions $a_1>0, a_2>0,a_3>0, a_1\left( a_2(a_3-a_6^2)-a_3a_5^2 \right)-a_4^2a_3a_2>0$. The condition $dF=0$ implies that $a_{14}=a_{24}=a_{34}=0.$
Then the condition $d\star F=0$ implies that $F$ is trivial.

\subsection{Lie algebra $\mathcal{A}_{3,9}\oplus\mathcal{A}_{1}$}
The structure of the Lie algebra is
\begin{equation*}
[e_1,e_2]=e_3,\quad [e_2,e_3]= e_1,\quad [e_3,e_1]= e_2. 
\end{equation*}
Up to automorphisms of the Lie algebra, a metric $g$ is given by
$$g=\left[\begin{array}{cccc}a_1 & 0 & 0 & a_4 \\0 & a_2 & 0 & a_5 \\0 & 0 & a_3& a_6 \\a_4 & a_5&a_6 & 1\end{array}\right],$$
with the conditions $a_1>0, a_2>0,a_3>0, a_1\left( a_2(a_3-a_6^2)-a_3a_5^2 \right)-a_4^2a_3a_2>0$. The condition $dF=0$ implies that $a_{14}=a_{24}=a_{34}=0.$ 
Then the only solution to $d\star F=0$ is the trivial solution.
\subsection{The abelian algebra}
Any metric can be reduced to the flat euclidean metric.

\section{Appendix}
We reproduce below the essential part of \textit{Maple} code used to solve the Einstein--Maxwell equations. We take here for example the Lie algebra $\mathcal{A}_{3,1}\oplus\mathcal{A}_{1}.$

{\footnotesize
\textit{\# To define the structure equations of the Lie algebra:}

\texttt{brac :=(x, y) $\rightarrow$ vector(n, [(x[3]*y[4]-x[4]*y[3]),0,0,0]);}

\textit{\# To define the metric:}

\texttt{G := Matrix([[a1, 0, 0, 0], [0, 1, 0, 0], [0, 0, 1, 0], [0, 0, 0, 1]])}

\textit{\# To define the coefficients of the Levi-Civita connection with respect to a metric $G$:}

\texttt{LeviCivita := (x, y, z) $\rightarrow$ (1/2)*evalm(innerprod(brac(x, y), G, z)+innerprod(brac(z, x), G, y)-innerprod(brac(y, z), G, x));}

\textit{\# To define the Levi-Civita connection of two vectors with respect to a $G$-orthonormal basis \{v[i]\}:}

\texttt{LC := (x, y) $\rightarrow$ evalm(LeviCivita(x, y, evalm(v[1]))*v[1]+LeviCivita(x, y, evalm(v[2]))*v[2]+LeviCivita(x, y, evalm(v[3]))*v[3]+LeviCivita(x, y, evalm(v[4]))*v[4]);}

\textit{\# To define the Riemannian tensor of the metric $G$:}

\texttt{Rc := (x, y, z) $\rightarrow$ -evalm(simplify(LC(x, LC(y, z))-LC(y, LC(x, z))-LC(brac(x, y), z)));}

\texttt{RiemC := (x, y, z, w) $\rightarrow$ simplify(innerprod(Rc(x, y, z), G, w));}

\textit{\# To define the Ricci tensor, the Riemannian scalar $R$ and the trace free part of the Ricci tensor $Ric_0$ of the metric $G$:}

\texttt{Ricci := (x, y)$\rightarrow$ evalm(simplify(RiemC(x, evalm(v[1]), y, evalm(v[1]))+RiemC(x, evalm(v[2]), y, evalm(v[2]))+RiemC(x, evalm(v[3]), y, evalm(v[3]))+RiemC(x, evalm(v[4]), y, evalm(v[4]))));}

\texttt{R := evalm(simplify(Ricci(evalm(v[1]), evalm(v[1]))+Ricci(evalm(v[2]), evalm(v[2]))+Ricci(evalm(v[3]), evalm(v[3]))+Ricci(evalm(v[4]), evalm(v[4]))));}

\texttt{Rico :=  (x, y)$\rightarrow$evalm(simplify(Ricci(x, y)-($\frac{1}{4}$)*R*innerprod(x, G, y)));}

\textit{\# To define the Hodge star of a $2$-form:}

\texttt{DGsetup([x1, x2, x3, x4]):}

\texttt{g := evalDG(a1*(dx1 \&t dx1)+(dx2 \&t dx2)+(dx3 \&t dx3)+(dx4 \&t dx4));}

\texttt{HodgeStar(g, a13*(dx1 \&w dx3)+a14*(dx1 \&w dx4)+a24(dx2 \&w dx4)+a23*(dx2 \&w dx3)+a34*(dx3 \&w dx4));}

\textit{\# To define the trace free part $[F\circ F]_0$ of a $2$-form $F$:}

\texttt{K := Transpose(Multiply(F, MatrixInverse(G)));}

\texttt{$F\circ F$ := simplify(Multiply(Transpose(K), F));}

\texttt{Tr($F\circ F$) := innerprod(evalm(v[1]), $F\circ F$, evalm(v[1]))+innerprod(evalm(v[2]), $F\circ F$, evalm(v[2]))+innerprod(evalm(v[3]), $F\circ F$, evalm(v[3]))+innerprod(evalm(v[4]), $F\circ F$, evalm(v[4]));}

\texttt{$[F\circ F]_0$ := simplify($F\circ F$-$\left(\frac{Tr(F\circ F)}{4}\right)$*G);}

\textit{\# To solve the Einstein--Maxwell equations:}

\texttt{solve({Rico(evalm(e[1]), evalm(e[1]))+ $[F\circ F]_0$[1, 1] = 0, Rico(evalm(e[1]), evalm(e[2]))+ $[F\circ F]_0$[1, 2] = 0, Rico(evalm(e[1]), evalm(e[3]))+ $[F\circ F]_0$[1, 3] = 0, Rico(evalm(e[1]), evalm(e[4]))+ $[F\circ F]_0$[1, 4] = 0, Rico(evalm(e[2]), evalm(e[2]))+ $[F\circ F]_0$[2, 2] = 0, Rico(evalm(e[2]), evalm(e[3]))+ $[F\circ F]_0$[2, 3] = 0, Rico(evalm(e[2]), evalm(e[4]))+ $[F\circ F]_0$[2, 4] = 0, Rico(evalm(e[3]), evalm(e[3]))+ $[F\circ F]_0$[3, 3] = 0, Rico(evalm(e[3]), evalm(e[4]))+ $[F\circ F]_0$[3, 4] = 0, Rico(evalm(e[4]), evalm(e[4]))+ $[F\circ F]_0$[4, 4] = 0}, {$a_1, a_{13}, a_{14}, a_{23}, a_{24}$});}
}

\bibliographystyle{abbrv}%\bibliographystyle{ieeetr}

\bibliography{biblio_E-M}
%\nocite{*}

\end{document}